\documentclass{ws-ijmpa}
\usepackage{amssymb, amsmath, amsfonts}
\usepackage[all]{xy}
\usepackage{bm}
\def\genfd{{\bm k}}

\long\def\nodo#1{{}}
\def\gg{\mathfrak{g}}

\def\antiu{S_{U(\gg)}} % antipode for U(\gg)
\def\PPartial#1{\frac{\partial}{\partial(\partial^{#1})}}

\def\End{\operatorname{End}}
\def\hx{\hat{x}}

\def\nxpoint{\refstepcounter{subsection}\makepoint{\thesubsection}}
\def\nxsubpoint{\refstepcounter{subsubsection}%
  \makepoint{\thesubsubsection}}
\def\refpoint#1{{\rm\textbf{\ref{#1}}}}
\def\makepoint#1{\medbreak\noindent{\bf #1. }}

\def\MR#1{} % we will not print MR numbers

\def\hpart{\hat\partial}
\begin{document}
%\begin{frontmatter}
\markboth{Zoran \v{S}koda} {Heisenberg double versus deformed
derivatives}
%%%%%%%%%%%%%%%%%%%%% Publisher's Area please ignore %%%%%%%%%%%%%%%
%
\catchline{}{}{}{}{}
%
%%%%%%%%%%%%%%%%%%%%%%%%%%%%%%%%%%%%%%%%%%%%%%%%%%%%%%%%%%%%%%%%%%%%
%\title{Heisenberg double versus deformed derivatives}
%\author{\sc Zoran \v{S}koda}
\title{HEISENBERG DOUBLE VERSUS DEFORMED DERIVATIVES}
\author{ZORAN \v{S}KODA}
\address{Theoretical Physics Division,
Institute Rudjer Bo\v{s}kovi\'{c}
\\ Bijeni\v{c}ka
cesta~54, P.O.Box 180, HR-10002 Zagreb, Croatia\\
{zskoda@irb.hr} }
\maketitle
\begin{history}
\received{Day Month Year} \revised{Day Month Year}
\end{history}

\begin{abstract} Two approaches to the tangent space of a
noncommutative space whose coordinate algebra is the enveloping
algebra of a Lie algebra are known: the Heisenberg double
construction and the approach via deformed derivatives, usually
defined by procedures involving orderings among noncommutative
coordinates or equivalently involving realizations via formal
differential operators. In an earlier work, we rephrased the
deformed derivative approach introducing certain smash product
algebra twisting a semicompleted Weyl algebra. We show here that
the Heisenberg double in the Lie algebra case, is isomorphic to
that product in a nontrivial way, involving a datum $\phi$
parametrizing the orderings or realizations in other approaches.
This way, we show that the two different formalisms, used by
different communities, for introducing the noncommutative phase
space for the Lie algebra type noncommutative spaces are
mathematically equivalent.\end{abstract}

%\vskip .1in
 \ccode{PACS 02.40.Gh; Mathematics
subject classification (2000): 16S30, 16S32, 16S35, 16Txx.}

%\vskip .03in
% 16S30 universal enveloping algebras
% 16S32 rings of differential operators
% 16S35 Twisted and skew group rings, crossed products
% 16Txx Hopf algebras, quantum groups and related topics
%Keywords: universal enveloping algebra, Hopf algebra, deformed
%derivative, Heisenberg double.

\vskip .1in

\section{Introduction}

\nxpoint Noncommutative algebras and noncommutative geometry may
play various roles in models of mathematical physics; for example
describing quantum symmetry algebras. A special case of interest
is when the noncommutative algebra is playing the role of the
space-time of the theory, and is interpreted as a small
deformation of the (commutative) 1-particle configuration space.
If one wants to proceed toward developing field theory on such a
space, it is beneficial to extend the deformation of the
configuration space to a deformation of full phase space
(symplectic manifold) of the theory. Deformed momentum space for
the noncommutative configuration space whose coordinate algebra is
the enveloping algebra of a finite-dimensional Lie algebra (also
called Lie algebra type noncommutative spaces) has been studied
recently in the mathematical physics
literature\cite{AC,DimGauge,HallidaySzabo}, mainly in special
cases, most notably variants of so-called $\kappa$-Minkowski
space\cite{BorPachol2,DimGauge,LukNowHeis,covKappadef,MS}.

\nxpoint ({\bf Deformed derivative approach})

\nxsubpoint (Notation) The algebras in the article are over a field $\genfd$ of
characteristic zero; both real and complex numbers appear in
applications of the present formalism. We fix a finite dimensional
Lie algebra $\gg$ with basis $\hat{x}_1,\ldots,\hat{x}_n$,
which are also the generators of
enveloping algebra $U(\gg)$; the corresponding commuting generators of the
symmetric algebra $S(\gg)$ will be denoted $x_1, \ldots,x_n$.

\nxsubpoint \label{pp:defDerAlg} Some authors (e.g.
\cite{DimGauge,HallidaySzabo,MS}) add to the (linear) enveloping
algebra generators, the corresponding ``deformed'' partial
derivatives $\partial^1,\ldots,\partial^n$, which are the dual
variables (in $\gg^*$), assumed to mutually commute. Then they
seek for the consistent commutation relations of the form

\begin{equation}\label{eq:commRel}
[\hat{x}_i, \partial^j] = \phi_i^j,\,\,\,\,\,\,\,i,j = 1,\ldots, n,
\end{equation}
where $\phi^i_j = \phi^i_j(\partial^1,\ldots,\partial^n)\in
\hat{S}(\gg^*)$ are formal power series and $\phi^i_j = \delta^i_j
+ $ higher order terms. By "consistency" they mean that one
quotients the free associative algebra product of $U(\gg)$ and of
the (commutative) formal power series ring in
$\partial^1,\ldots,\partial^n$ (the latter is isomorphic to the
completion of polynomial ring in the dual variables
$\hat{S}(\gg^*)$) by the commutation relations (\ref{eq:commRel})
and the restriction of the quotient map to each of the parts,
$U(\gg)$ and $\hat{S}(\gg^*)$ separately, has a zero kernel. The
resulting quotient algebra generated by
$\hat{x}_1,\ldots,\hat{x}_m,\partial^1,\ldots,\partial^n$ will be
referred to as the {\bf phase space algebra with the
$\phi$-deformed derivatives} (of the noncommutative algebra
$U(\gg)$). It follows from the Jacobi identities\cite{scopr} that
this nondegeneracy condition for the matrix $(\phi^i_j)_{i,j
=1,\ldots,n}$ can be expressed by requiring that $(\phi^i_j)$
provides a solution to the system
\begin{equation}\label{eq:phiderphi}
\phi^l_j \PPartial{l}(\phi^k_i) - \phi^l_i \PPartial{l}(\phi^k_j)
= C^s_{ij}\phi^k_s.
\end{equation}
of formal differential equations (summation on repeated indices
understood). Moreover a solution exists for all $\gg$ (such a {\it
universal} solution is exhibited in \cite{ldWeyl}), but the
solution is not unique; moreover, we have shown in \cite{scopr}
that if we require $\phi^i_j = \delta^i_j +$ higher order terms,
then the choice of such a solution $\phi^i_j$ is equivalent to any
among some other data of interest (some of the equivalences known
before):

- ``ordering prescription''\cite{HallidaySzabo,MS}

- realization of the enveloping algebra in a semi-completed Weyl
algebra\cite{ldWeyl,scopr} of the form $\hat{x}^\phi_i = \sum_j
x_j \phi^j_i$;

- a homomorphism of Lie algebras $\bm\phi:\gg\to
Der(\hat{S}(\gg^*))$ (then $\phi^j_i =
\bm\phi(-\hat{x}_i)(\partial^j)$). It extends to a Hopf action
also denoted $\bm\phi:U(\gg)\to End(\hat{S}(\gg^*))$.

- a deformed Leibniz rule showing how $\partial^i$ acts on a
product $\hat{u}\hat{v}$ of $\hat{u},\hat{v}\in U(\gg)$);

- (topological) coproduct $\Delta: \hat{S(\gg)^*}\to
\hat{S}(\gg^*)\hat\otimes\hat{S}(\gg^*)$;

- prescription for multiplying certain formal exponentials of a
noncommutative argument\cite{exp,kappaind,HallidaySzabo} (not
shown in \cite{scopr});

- a choice of the star products (belonging to a specific class of
star products);

- a coalgebra isomorphism $\xi: S(\gg)\to U(\gg)$ such that
$\xi|_\gg = \mathrm{id}_\gg$;

For the purpose of the proofs we sketch below some of the relations
among the above data, for more see \cite{scopr} and Sec. \ref{sec:defDer}.

\nxpoint \label{smash} ({\bf Hopf actions and smash products})
Recall that a left action $\triangleright :H\otimes A\to A$ of a
Hopf algebra $H$ on an algebra $A$ is a {\bf Hopf action} if it is
satisfying the condition $h\triangleright (a\cdot b) = \sum
(h_{(1)}\triangleright a)\cdot(h_{(2)}\triangleright b)$, where we
used the Sweedler notation $\Delta(h)=\sum h_{(1)}\otimes
h_{(2)}$; we also say that $A$ is a left $H$-module algebra. In
that case, one defines the smash product algebra (or crossed
product) $A\sharp H$ as the tensor product $A\otimes H$ with the
associative multiplication given by
$$(a\otimes h)(b\otimes g) =
\sum (a h_{(1)}\triangleright b)\otimes (h_{(2)} g).$$

\nxpoint ({\bf Heisenberg double}) The input for the Heisenberg
double\cite{KaprHeis,Lu,Semikhatov} construction is a pair of Hopf
algebras $H,H'$ in a bilinear pairing $\langle,\rangle:H\otimes
H'\to\genfd$ which is Hopf, i.e. with the product on pairings on
the tensor square, the coproduct and the product are dual in the
sense $\langle \Delta_H(a),b\otimes c\rangle = \langle a, b\cdot
c\rangle$, $\langle a\otimes a',\Delta_{H'}b\rangle = \langle
a\cdot a', b\rangle$ and similarly for the unit and counit. In our
case $H=U(\gg)$ and the role of $H'$ is played by the algebraic
linear dual $U^*(\gg)=\mathrm{Hom}_\genfd(U(\gg),\genfd)$ which is
a {\it topological Hopf algebra}, i.e. the coproduct of the
generators may result in infinitely many summands from the tensor
square, amounting to the need for some completion of $H\otimes
H'$. Similar to the Drinfel'd double, Heisenberg double is the
algebra whose underlying space is (a completion of) $H\otimes H'$,
but unlike Drinfel'd double it does not have a Hopf algebra
structure itself. One defines the {\bf coregular action} of $H'$
on $H$ given by $h'\triangleright h = \sum h_{(1)}\langle h_{(2)},
h'\rangle $ where $\Delta_H(h)=\sum h_{(1)}\otimes h_{(2)}$; as we
required that the pairing is Hopf pairing, this action of $H'$ on
$H$ is automatically a Hopf action (cf. \refpoint{smash}), hence
we can form the corresponding smash product algebra $H\sharp H'$,
the Heisenberg double of $H$ (better, of the data
$(H,H',\langle,\rangle)$).

\nxpoint ({\bf Sketch of the proof of the main result})
\label{p:sketch}

We want to exhibit the isomorphism between the phase space algebra
with $\phi$-deformed derivatives, and the Heisenberg double of
$U(\gg)$. This comprises four steps/isomorphisms, the first two of
which were effectively done in our earlier work\cite{scopr}, and
the remaining step is the focus of this paper.

I By the definition, the phase space algebra with the
$\phi$-deformed derivatives is generated by
$\hat{x}_1,\ldots,\hat{x}_n$ in $U(\gg)$ and the mutually
commuting formal power series in $\partial^1,\ldots,\partial^n$
with commutation relation~(\ref{eq:phiderphi}). In \cite{scopr},
we have shown that it is isomorphic to the smash product
$U(\gg)\sharp_\phi \hat{S}(\gg^*)$.

II By~\cite{scopr}, $\phi$ induces an isomorphism of coalgebras
denoted $\xi_\phi : S(\gg)\to U(\gg)$. The transpose of $\xi_\phi$
is an isomorphism of topological algebras $\xi^T_\phi : S(\gg)^*
\to U(\gg)^*$ which one composes with the isomorphism
$\hat{S}(\gg^*)\cong S(\gg)^*$. Therefore the algebra isomorphism
$U(\gg)\sharp_\phi \hat{S}(\gg^*)\cong U(\gg)\sharp U(\gg)^*$
where the action used for the smash product also transfers to the
right hand side. This isomorphism is also exhibited in
\cite{scopr}. Notice that the smash product on the right hand side
is {\it not} yet the Heisenberg double as the action used is the
action of $U(\gg)$ on $U(\gg)^*$ and not the conversely.

IIB Note that the isomorphism $\hat{S}(\gg^*)\cong U(\gg)^*$
obtained via $\xi^T_\phi$ (cf. II) and the identification
$\hat{S}(\gg^*)\cong S(\gg)^*$ induces also a nondegenerate Hopf
pairing between $\hat{S}(\gg^*)$ and $U(\gg)$. For this pairing we
find several descriptions~(\ref{eq:pairing1}),(\ref{eq:pairing2})
which are used below to describe the Heisenberg double of
$U(\gg)$.

III The smash product algebra depends on the action used in its
definition. The smash product $U(\gg)\sharp_\phi \hat{S}(\gg^*)$
in II is derived from the action $\triangleright$ of the Hopf
algebra $U(\gg)$ on $\hat{S}(\gg^*)$. We relate this action with
the ``black'' action $\blacktriangleright$ (see
\ref{p:blackAction}) of the topological algebra
$\hat{S}(\gg^*)\cong U(\gg)^*$, and show that the two resulting
smash products (one from action of $U(\gg)$ on $S(\gg)$ and
another from the black action of $U(\gg^*)$ on $U(\gg)$) are
isomorphic as abstract algebras. The action $\blacktriangleright$
is a Hopf action with respect to the topological coproduct on
$U(\gg)^*$ or equivalently the $\phi$-deformed coproduct on
$\hat{S}(\gg^*)$.

IV We show in \refpoint{p:main} that the black action
$\blacktriangleright$ is precisely the coregular action, i.e. the
unique action satisfying $P\blacktriangleright \hat{u} = \sum
\hat{u}_{(1)}\,\langle \hat{u}_{(2)}, P\rangle_{\bm\phi}$ where
$\langle , \rangle_{\bm\phi}$ is the Hopf pairing with the dual
topological Hopf algebra (with the dual represented in a specific
way). The coregular action is used in the definition of the
Heisenberg double, completing the identification
$U(\gg)\sharp_\phi \hat{S}(\gg^*)\cong U(\gg)\sharp U(\gg)^*$
where for the smash products, on the left hand side one uses the
$U(\gg)$ action by $\triangleright$, and on the right hand side
the coregular $U(\gg)^*$-action.

\section{More on deformed derivatives}

\label{sec:defDer}
More familiarity with the structure involved in the method of
deformed derivatives is needed later to exhibit its relation
to the Heisenberg double. For the users of our results we also
sketch the connection to star products.

\nxpoint ({\bf Star product perspective}) Lie algebra type
noncommutative spaces are simply the deformation quantizations of
the linear Poisson structure; given structure constants $C^k_{ij}$
linear in a deformation parameter the enveloping algebras of the
Lie algebra $\gg$ given in a base by $[\hat{x}_i,\hat{x}_j]=
C^k_{ij}\hat{x}_k$ is viewed as a deformation of the polynomial
(symmetric) algebra $S(\gg)$ generated by commuting
$x_1,\ldots,x_n$. Given any linear isomorphism $\xi:
S(\gg)\stackrel\cong\longrightarrow U(\gg)$ we transfer the
noncommutative product on $U(\gg)$ to a $\star$-product on
$S(\gg)$, defined by $f\star g = \xi^{-1}(\xi(f)\cdot\xi(g))$.
There are many isomorphisms which may play role of $\xi$, but in
order to introduce either the $\phi$-deformed derivatives like
in~\cite{AC,HallidaySzabo,MS,scopr}, or to make the correspondence
with the Heisenberg double construction, we need to restrict to
$\xi$ which are {\em coalgebra isomorphisms}; we also require a
``small deformation condition'' that $\xi$ is the identity on the
constant and linear parts, i.e. on $\genfd\oplus\gg\subset
S(\gg)$. Our restriction to coalgebra isomorphisms, singles out a
distinguished class of star products quantizing the linear Poisson
structure. Kathotia\cite{Kathotia} compares the Kontsevich star
product\cite{KontsDefPoiss} for linear Poisson structures to the
PBW-product which corresponds to the case where $\xi$ is the
standard symmetrization (coexponential) map (cf. \cite{ldWeyl},
especially Chapter 10); Kontsevich star product is not in our
class, although it is equivalent to the PBW product, which is in
our class.

\nxpoint ({\bf Some connections between the basic data}) Coalgebra
isomorphism $\xi:S(\gg)\to U(\gg)$ induces a transpose map
$\xi^T:U^*(\gg)\to S^*(\gg)$, which is consequently an algebra
isomorphism. There is an isomorphism $S^*(\gg)\cong\hat{S}(\gg^*)$
where $\hat{S}(\gg^*)$ denotes a completed symmetric algebra on
the dual; the isomorphism depends on a normalization of a pairing
between $\hat{S}(\gg^*)$ and $S(\gg)$ (cf.~\cite{ldWeyl}, {\bf
10.4, 10.5}). the functionals in $S^*(\gg)\cong \hat{S}(\gg^*)$
can be identified with the infinite order differential operators
with constant coefficients: a differential operator applied to a
polynomial in $S(\gg)$ and then evaluated at $0$, defines a
differential operator. If the dual generators of $\gg^*\subset
\hat{S}(\gg^*)$ corresponding to the basis $x_1,\ldots,x_n$ are
denoted as the partial derivatives $\partial^i$, this rule and
identification explains the choice of normalization
in~\cite{ldWeyl}, Sec. {\bf 10}. The topological coproduct on
$U^*(\gg)$ which is the algebraic transpose to the product on
$U(\gg)$, is (for $\xi$ being the symmetrization map) written as a
formal differential operators in $\hat{S}(\gg^*)$
in~\cite{ReshQLieBialg}, where the generalizations for Lie
bialgebras are considered. In~\cite{scopr} we have shown that this
deformed coproduct is the same as a coproduct obtained by using
Leibniz rules defined in terms of the deformed commutation
relations; and in the case of symmetric ordering we have
exhibited\cite{scopr} a Feynman-like diagram expansion summing to
what is essentially a Fourier-transformed form of the BCH series.

\nxpoint \label{phidata} As shown in detail in \cite{scopr}, the
coalgebra isomorphism $\xi:S(\gg)\to U(\gg)$ tautological on
$\genfd\oplus\gg$ as above, is equivalent to any of several other
data listed in the introduction, including the $\phi$-data
described as follows. The star product $x_i \star f$ is always of
the form $\sum_j x_j \phi^i_j(\partial^1,\ldots,\partial^n)(f)$
where $(\phi^i_j)_{i,j = 1,\ldots,n}$ is a matrix of elements in
$\hat{S}(\gg^*)$ (formal power series in dual variables
$\partial^1,\ldots,\partial^n$) satisfying a formal set of
differential equations (\cite{ldWeyl} ch. 4) equivalent to the
statement that the formula
$\bm\phi(-\hat{x}_i)(\partial^j)=\phi_i^j$ defines a Lie algebra
morphism $\bm\phi:\gg\to \mathrm{Der}(\hat{S}(\gg^*))$.

The correspondence $\hat{x}_i\mapsto \hat{x}_j^\phi = \sum_j x_j
\phi^j_i$ extends to an injective morphism of associative algebras
$()^\phi : U(\gg)\to \hat{A}_{n,\genfd}$ where
$\hat{A}_{n,\genfd}$ is the Weyl algebra of differential operators
with polynomial coefficients, completed by the degree of the
differential operator (hence we allow formal power series in
$\partial^i$-s but not in $x_j$-s). This (semi)completed Weyl
algebra has the standard Fock representation on $S(\gg)$. The Lie
algebra homomorphism $\bm\phi$ extends multiplicatively to a
unique homomorphism $U(\gg)\to\mathrm{End}(\hat{S}(\gg^*))$ (also
denoted $\bm\phi$), which is a Hopf action (cf. \refpoint{smash}).
Thus we can form a smash product algebra $A_{\gg,\phi} =
U(\gg)\sharp\hat{S}(\gg^*)$, the semicompleted $n$-th Weyl algebra
($\hat{A}_{n,\genfd}$ is the special case of this construction for
an abelian Lie algebra). The rule
$\hat{x}_i\mapsto\hat{x}^\phi_i$, $\partial^j\to\partial^j$
extends to a unique homomorphism
$A_{\gg,\phi}\to\hat{A}_{n,\genfd}$; one easily shows that it is
an isomorphism.

\nxpoint \label{p:blackAction}({\bf The action later used for
Heisenberg double}) Not only $U(\gg)$ acts by Hopf action on
$\hat{S}(\gg^*)$ (this action was used in the construction of
$A_{\gg,\bm\phi}$), but also conversely $\hat{S}(\gg^*)$ as a
topological Hopf algebra acts on $U(\gg)$. The {\em latter} action
$\blacktriangleright$ is in the Main Theorem below identified with
the smash product action of the Heisenberg double! To define the
latter action, $U(\gg)$ is embedded as a subalgebra $U(\gg)\sharp
\genfd\hookrightarrow A_{\gg,\phi}$; and similarly for
$\hat{S}(\gg^*)$. The action $a\otimes u\mapsto
a\blacktriangleright \hat{u}$, $A_{\gg,\bm\phi}\otimes U(\gg)\to
U(\gg)$ is defined by multiplying within $A_{\gg,\phi}$ and then
projecting by evaluating the second tensor factor in
$A_{\gg,\phi}=U(\gg)\sharp\hat{S}(\gg^*)$ (as a differential
operator) at $1$. Thus $U(\gg)$ is an $A_{\gg,\bm\phi}$-module,
the deformed Fock space where $1_{U(\gg)}$ is the $\phi$-deformed
vacuum. It can be shown\cite{scopr} that the coalgebra isomorphism
$\xi:S(\gg)\to U(\gg)$ can be computed by composing
$S(\gg)\hookrightarrow \hat{A}_{n,\genfd}\cong
A_{\gg,\bm\phi}\stackrel{\blacktriangleright
1_{U(\gg)}}\longrightarrow U(\gg)$.

\nxpoint ({\bf Deformed coproduct}) If we define, for
$P\in\hat{S}(\gg^*)\hookrightarrow A_{\gg,\bm\phi}$, the linear
operator $\hat{P}:U(\gg)\to U(\gg)$ by
$\hat{P}(\hat{u})=P\blacktriangleright\hat{u}$ or (equivalently,
according to~\cite{scopr}) by $\hat{P}(\xi(f)) = \xi(P(f))$, then
the Leibniz rule holds:
$\sum\hat{P}_{(1)}(\hat{u})\cdot_{U(\gg)}\hat{P}_{(2)}(\hat{v}) =
P(\hat{u}\cdot_{U(\gg)}\hat{v})$ for a unique ($\phi$-dependent)
deformed coproduct $P\mapsto \Delta(P) = \sum P_{(1)}\otimes
P_{(2)}$ on $\hat{S}(\gg^*)$ (with the tensor product allowing
infinitely many terms), cf. \refpoint{p:Leibnizdef}.

\section{Relating Heisenberg double to the $\phi$-deformed derivatives}

\nxpoint {\bf Lemma.} {\it The following nonsymmetric formula for
$\Delta(\partial^\mu)$ holds}:

\begin{equation}\label{eq:nonsym}
\Delta(\partial^\mu) = 1\otimes \partial^\mu+\partial^\alpha
\otimes [\partial^\mu,\hat{x}_\alpha] +
\frac{1}{2!}\partial^{\alpha_1}\partial^{\alpha_2}\otimes
[[\partial^\mu,\hat{x}_{\alpha_1}],\hat{x}_{\alpha_2}] +\ldots
\end{equation}
The sum has only finitely many terms when applied to an element in
$U(\gg)\otimes U(\gg)$. Proof is by induction, see~\cite{scopr}.

\nxpoint {\bf Lemma.} {\it If $\hat{a} = \sum_{\alpha = 1}^n
a^\alpha \hat{x}_\alpha$ and $\hat{f}\in U(\gg)$ then}
\begin{equation}\label{eq:dxn}
\hpart^\mu (\hat{a}^p \hat{f}) = \sum_{k = 0}^{p-1} {p\choose k}
a^{\alpha_1} a^{\alpha_2}\cdots a^{\alpha_k} \hat{a}^{p-k}
[[\hpart^\mu,\hx_{\alpha_1}],\ldots,\hx_{\alpha_k}] (\hat{f})
\end{equation}
{\it Proof.} This is a tautology for $p = 0$. Suppose it holds for
all $p$ up to some $p_0$, and for all $\hat{f}$. Then set $\hat{g}
= \hat{a}\hat{f} = a^\alpha \hat{x}_\alpha$. Then
$\hpart^\mu(\hat{a}^{p_0 +1}\hat{f})
=\hpart^\mu(\hat{a}^{p_0}\hat{g})$ and we can apply~(\ref{eq:dxn})
to $\hpart^\mu(\hat{a}^{p_0}\hat{g})$. Now
$$\begin{array}{lcl}
[[[\hpart^\mu,\hx_{\alpha_1}],\ldots],\hx_{\alpha_k}] (\hat{g})
&=&
a^{\alpha_k}[[[\hpart^\mu,\hx_{\alpha_1}],\ldots],\hx_{\alpha_k}]
(\hx_{\alpha_{k+1}}\hat{g})\\
&=&\hat{a}[[[\hpart^\mu,\hx_{\alpha_1}],\ldots],\hx_{\alpha_k}]
(\hat{f}) \,+ \\&&\,\,\,\,\,\,\,+ \,\,a^{\alpha_{k+1}}
[[[[\hpart^\mu,\hx_{\alpha_1}],\ldots],\hx_{\alpha_k}],\hx_{\alpha_{k+1}}]
(\hat{f}).
\end{array}$$
Collecting the terms and the Pascal triangle identity complete the
induction step.

\nxsubpoint {\bf Remark.} It is interesting that this lemma was
needed and proved in~\cite{scopr} related to certain Feynman
diagram type expansion calculation leading to an exact summation
result, whereas it will be seen here rather as a step toward and a
special case of a formula showing the condition that certain
secondary action in the $\phi$-deformed derivatives picture (the
black action) is precisely the coregular action needed to define
the Heisenberg double.

\nxpoint {\bf Theorem.} {\it Given a left Hopf action
$\bm\phi:U(\gg)\to\End(\hat{S}(\gg^*))$, with
$\bm\phi(-\hx_i)(\partial^j)=\delta^i_j + O(\partial)$, there is a
Hopf pairing $\langle , \rangle_{\bm\phi} : U(\gg)\otimes
\hat{S}(\gg^*)\to\genfd$ given by
\begin{equation}\label{eq:pairing1}
\langle \hat{u},P\rangle_{\bm\phi} =
\bm\phi(\antiu\hat{u})(P)|0\rangle \equiv
\bm\phi(\antiu\hat{u})(P) (1_{S(\gg)})
\end{equation}
where $\hat{u}\in U(\gg)$, $P\in\hat{S}(\gg^*)$, and $\antiu$ is
the antipode antiautomorphism of $U(\gg)$, and where
$\hat{S}(\gg^*)$ is considered a topological Hopf algebra with
respect to the $\bm\phi$-deformed coproduct. }

{\it Proof.} Clearly the pairing is well defined; the antipode
comes because we use {\it left} Hopf actions. The product of
diferential operators {\it with constant coeficients} evaluated at
$1$ equals the product of their evaluations at $1$. Therefore the
fact that $\bm\phi$ is Hopf action implies $\langle
\hat{u},PQ\rangle_{\bm\phi} = \langle\Delta\hat{u},P\otimes
Q\rangle_{\bm\phi}$. It is less obvious to verify the other
duality: of $\bm\phi$-deformed coproduct and the multiplication on
$U(\gg)$. It is sufficient to show that one has
\begin{equation}\label{eq:stepduality}
\langle\hat{x}_\alpha\hat{u},\partial^\mu\rangle_{\bm\phi} =
\langle\hat{x}_\alpha\otimes\hat{u},\Delta\partial^\mu\rangle_{\bm\phi}.
\end{equation}
for all $\alpha$ and all $\hat{u}$ in $U(\gg)$. Indeed, extending
to $\prod_{i=1}^k x_{\alpha_i}\hat{u}$ for all
$(\alpha_1,\ldots,\alpha_k)$ can be done by induction on $k$,
using the coassociativity of the coproduct and associativity of
the product. Once it is true for any product $\hat{v}\hat{u}$ in
the left argument, it is an easy general nonsense, using the
already known duality for $\Delta_{U(\gg)}$, to extend the
property to products of $\partial$-s by induction using the
following calculation for the induction step
$$\begin{array}{lcl}
\langle\hat{v}\hat{u},P_1 P_2\rangle_{\bm\phi} & = &\langle
\sum\hat{v}_{(1)} \hat{u}_{(1)}\otimes \hat{v}_{(2)}\hat{u}_{(2)},
P_1\otimes P_2\rangle_{\bm\phi}
\\&=& \sum\langle\hat{v}_{(1)}\otimes\hat{u}_{(1)},\Delta(P_1)\rangle_{\bm\phi}
\langle\hat{v}_{(2)}\otimes\hat{u}_{(2)},\Delta(P_2)\rangle_{\bm\phi}
\\&=&\sum\langle\hat{v}_{(1)}\otimes\hat{u}_{(1)}\otimes
\hat{v}_{(2)}\otimes\hat{u}_{(2)},
\Delta(P_1)\otimes\Delta(P_2)\rangle_{\bm\phi}
\\&=&\sum\langle\hat{v}\otimes\hat{u},
\Delta(P_1 P_2)\rangle_{\bm\phi}
\end{array}$$
Let us now calculate~(\ref{eq:stepduality}) using the nonsymmetric
formula~(\ref{eq:nonsym}) for the $\bm\phi$-coproduct. All terms
readily give zero in first factor unless the first factor is
degree $1$ in $\partial$-s. Thus we effectively need to show
$$
\sum_\beta \langle x_\alpha, \partial^\beta\rangle_{\bm\phi}
\otimes
\langle\hat{u},[\partial^\mu,\hat{x}_\beta]\rangle_{\bm\phi} =
\langle x_\alpha\hat{u},\partial^\mu\rangle_{\bm\phi}.
$$
The left-hand side is $\sum_\beta
\bm\phi(-\hat{x}_\alpha)(\partial^\beta)
\bm\phi(\antiu\hat{u}^{\mathrm{op}})(\bm\phi(-\hat{x}_\beta)(\partial^\mu))|0\rangle
=$\newline $=\sum_\beta
\bm\phi(-\hat{x}_\alpha)(\partial^\beta)|0\rangle
\bm\phi(\antiu(\hat{u}^{\mathrm
{op}})\hat{x}_\beta)(\partial^\mu))|0\rangle$ and
$\bm\phi(-\hat{x}_\alpha)(\partial^\beta)|0\rangle =
\delta^\beta_\alpha$ by the assumption on $\bm\phi$. Finally, the
contraction with the Kronecker delta gives
$\bm\phi(\antiu(\hat{x}_\alpha\hat{u})^{\rm
op})(\partial^\mu)|0\rangle$.

\nxpoint {\bf Proposition.} {\it If $\xi=\xi_{\bm\phi}:S(\gg)\to
U(\gg)$ is the coalgebra isomorphism correspoding to $\bm\phi$ and
$\xi^T:U(\gg)^*\to S(\gg)^*\cong \hat{S}(\gg^*)$ its transpose,
then the pairing may be described alternatively by
\begin{equation}\label{eq:pairing2}
\langle \hat{u},P\rangle_{\bm\phi} =  (\xi^T)^{-1}(P)(\hat{u}) =
P(\xi^{-1}_{\bm\phi}(\hat{u})) =
\epsilon_{S(\gg)}(P(\hat{u}^{\bm\phi}|0\rangle)) =
\epsilon_{S(\gg)}((\hat{P}(\hat{u}))^\phi|0\rangle),
\end{equation}
where $P (\xi^{-1}_{\bm\phi}(\hat{u}))$ is the evaluation of
$P\in\hat{S}(\gg^*)$ on $\xi^{-1}_{\bm\phi}(\hat{u})\in S(\gg)$
via the pairing. } \vskip .07in

{\it Proof.} We show $\langle \hat{u},P\rangle_{\bm\phi} =
\epsilon_{S(\gg)}(P(\hat{u}^{\bm\phi}|0\rangle))$. By the previous
arguments, it is enough to show that this alternative formula
gives the same (and, in particular, Hopf) pairing
as~(\ref{eq:pairing1}) when $P=\partial^\mu$. This is evident when
$\hat{u} = \hx_\nu$ for some $\nu$. Now suppose by induction
that~(\ref{eq:pairing2}) holds for $\hat{u}$. Then
$$\begin{array}{lcl}
\epsilon\partial^\mu(\hx_\lambda^\phi\hat{u}^\phi|0\rangle) &=&
\epsilon\partial^\mu x_\alpha\phi^\alpha_\lambda \hat{u}^\phi
|0\rangle
\\ &=& \epsilon x_\alpha \partial^\mu \phi^\alpha_\lambda \hat{u}^\phi |0\rangle
+ \epsilon \phi^\mu_\lambda \hat{u}^\phi |0\rangle \\&=& 0 +
\epsilon\bm\phi(\antiu\hat{u})([\hat\partial^\mu,\hx_\lambda]) \\
&=& \epsilon\bm\phi(\antiu(\hx_\lambda\hat{u}))(\partial^\mu),
\end{array}$$
hence it holds for $\hx_\lambda\hat{u}$.

The other equalities in~(\ref{eq:pairing2}) are direct:
$(\xi^T)^{-1}(P)(\hat{u}) = P(\xi^{-1}_{\bm\phi}(\hat{u}))$ by the
definition of the transpose operator $\xi^T$; then
$P(\xi^{-1}_{\bm\phi}(\hat{u})) = \epsilon
(P(\hat{u}^{\bm\phi}|0\rangle))$ and $\epsilon
(P(\hat{u}^{\bm\phi}|0\rangle)) = \epsilon
((\hat{P}(\hat{u}))^\phi|0\rangle)$ by the basic identities
$\epsilon (\hat{w}|0\rangle) = \xi^{-1}(\hat{w})$ and
$\hat{P}\circ\xi = \xi\circ P$.

\nxpoint \label{p:main} {\bf Main Theorem.} {\it The
$(\gg,\bm\phi)$-twisted Weyl algebra $A_{\gg,\bm\phi}$ is
isomorphic to the Heisenberg double of the Hopf algebra $U(\gg)$
where the dual topological Hopf algebra is $\hat{S}(\gg^*)$ with
respect to the $\bm\phi$-deformed coproduct, and with respect to
the Hopf pairing given by ~(\ref{eq:pairing1}) or, equivalently,
~(\ref{eq:pairing2}). In other words, the left action
$\blacktriangleright$ used for the second smash product structure
satisfies (and is determined by) the formula
$$
P\blacktriangleright \hat{u} = \sum \langle \hat{u}_{(2)},
P\rangle_{\bm\phi}\, \hat{u}_{(1)}
$$
for all  $\hat{u}\in U(\gg)$ and $P\in\hat{S}(\gg^*)$.

Consequently, the phase space algebra with the $\phi$-deformed
derivatives~(\refpoint{pp:defDerAlg}) is isomorphic to the
Heisenberg double (and the isomorphism nontrivially depends on
$\phi$). }

 {\it Proof.} If the identity holds for $P=P_1$ and $P=P_2$ then
$$\begin{array}{lcl}
P_1 P_2\triangleright\hat{u}&=& P_1\triangleright\left(\sum\langle
\hat{u}_{(2)},P_2\rangle_{\bm\phi}\, \hat{u}_{(1)}\right)\\
 &=& \sum\langle \hat{u}_{(3)},P_2\rangle_{\bm\phi}\, \langle
\hat{u}_{(2)},P_1\rangle_{\bm\phi}\, \hat{u}_{(1)}\\
&=& \sum\langle \hat{u}_{(2)},P_1 P_2\rangle_{\bm\phi}\,
\hat{u}_{(1)}
\end{array}$$
hence it holds for $P=P_1 P_2$. For $P=1$ it holds trivially,
hence it is sufficient to check for $P=\partial^\mu$ and use
induction. The identity is linear in $\hat{u}\in U(\gg)$, so it is
sufficient to prove it for all $\hat{u}$ of the form $\hat{u} =
\hat{a}^p = (\sum_{\alpha = 1}^n a^\alpha \hat{x}_\alpha)^p$,
$p\geq 0$ where $\hat{a} = \sum_\alpha a^\alpha \hat{x}_\alpha$ is
arbitrary. In that case, $\Delta(\hat{u}) = \sum_{k=0}^{p}
{p\choose k} \hat{a}^{p-k}\otimes \hat{a}^{k}$ and we need to show
$$
\partial^\mu\blacktriangleright \hat{a}^p
= \hat{\partial}^\mu (\hat{a}^p) = \sum_{k=0}^p {p \choose k}
\langle \hat{a}^k, \partial^\mu\rangle_{\bm\phi}\, \hat{a}^{n-k}
$$
but $\langle \hat{a}^k, \partial^\mu\rangle_{\bm\phi}$ is
by~(\ref{eq:pairing2}) equal to
$$\bm\phi(S_{U(\gg)}(\hat{a}^k))(\partial^\mu)
= (-1)^k\bm\phi(\hat{a}^k)(\partial^\mu) = (-1)^k
[\ldots[[\partial^\mu,\hat{a}],\hat{a}],\ldots,\hat{a}],$$ what by
linearity reduces to~(\ref{eq:dxn}) for the case $f=1$. (This
shows IV in \ref{p:sketch} i.e. that the black action
$\blacktriangleright$ is identifiable with the coregular action
under the isomorphism $\hat{S}(\gg^*)\cong U(\gg)^*$).

While the vector spaces of the two smash products
($U(\gg)\sharp_\phi \hat{S}(\gg)$ and the Heisenberg double) are
isomorphic by the definition (they are simply the tensor products
with the same factors), we need to show that the multiplication is
the same; for this we need to commute the tensor factors. One can
easily compute that if $[\partial,\hat{x}] = Q\in \hat{S}(\gg^*)$,
then also in the Heisenberg double $[\hat\partial,\hat{x}] =
\hat{Q}$ for $\partial\in\gg^*$ and $\hat{x}\in\gg\hookrightarrow
U(\gg)$. Therefore for the generators, the commutation relations
in the two smash products agree (this shows III in
\ref{p:sketch}), hence the isomorphism of $A_{\gg,\bm\phi}$ and
the smash product given by the black action, hence the Heisenberg
double.

The final sentence in the theorem now follows by I in
\refpoint{p:sketch}, namely we know from our earlier
work\cite{scopr} that the phase space algebra with the deformed
derivatives is isomorphic to the smash product $A_{\gg,\bm\phi}$.
Step II in \refpoint{p:sketch} shown in \cite{scopr} is used all
along in the construction. Notice that the heart of this paper is
performing the step IV from \refpoint{p:sketch}; once we have done
it, we have recapitulated earlier prepared steps for I, II and
III.

\section{Final remarks.}

\nxpoint Though the $\phi$-deformed derivatives are not present
there, the Reshetikhin's article\cite{ReshQLieBialg} has
implicitly much of the structure from this paper (including issues
on dualization of coproducts) implicitly present, including the
quantum deformations of enveloping algebras and more general
bialgebras.

\nxpoint \label{p:Leibnizdef} The fact that the Leibniz rule for
the action of $\hat{S}(\gg^*)$ on $U(\gg)$ (for any $\gg$ and
$\bm\phi$) gives a well-defined coassociative map into the tensor
product is not obvious in the deformed derivative
picture\cite{AC,DimGauge,scopr,MS}; namely it is {\it a priori}
undefined up to a kernel of the multiplication map (add an element
in the kernel and the Leibniz rule does not change). But now the
Hopf action is well-defined within the Heisenberg double
construction and the Heisenberg double as an algebra is identified
with $A_{\gg,\bm\phi}$ where the deformed Leibniz rule was
originally defined. Heisenberg double provides an invariant
picture, giving simple "dual" interpretation to the deformed
coproduct, while the approach via the $\phi$-deformed derivatives
and commutators is useful for calculation, as it is exhibited in
the physics literature before.

\nxpoint The differential forms and exterior derivative can also
be extended to the same setup\cite{exterior,meljDiff}.

\vskip .1in {\footnotesize {\bf Acknowledgments} Main results were
obtained at the Institute Rudjer Bo\v{s}kovi\'c, Zagreb; a part of
the article has been written at Max Planck Institute for
Mathematics, Bonn, whom I thank for excellent working conditions;
travel has been partly supported by Croatian (MZO\v{S})-German
(DAAD) bilateral project organized together with Urs Schreiber.}

\footnotesize{

}

\begin{thebibliography}{99}
\bibitem{AC} G. Amelino-Camelia,
M. Arzano, Coproduct and star product in field theories on
Lie-algebra non-commutative space-times, {\em Phys. Rev.}
D65:084044 (2002) {\tt hep-th/0105120}.

\bibitem{BorPachol2}
A. Borowiec, A. Pacho\l, $\kappa$-Minkowski spacetimes and DSR
algebras: fresh look and old problems, {\em SIGMA} 6:086 (2010),
{\tt arxiv/1005.4429}

\bibitem{DimGauge}
M. Dimitrijevi\'c, F. Meyer, L. M\"oller, J. Wess , Gauge theories
on the $\kappa$-Minkowski spacetime, {\em Eur. Phys.J.} C36 (2004)
117--126; {\tt hep-th/0310116}.

\bibitem{ldWeyl}
N. Durov, S. Meljanac, A. Samsarov, Z. \v{S}koda, A universal
formula for representing Lie algebra generators as formal power
series with coefficients in the Weyl algebra, {\em J. Algebra}
309:1, 318--359 (2007) (math.RT/0604096)

\bibitem{HallidaySzabo}
S. Halliday, R. J. Szabo, Noncommutative field theory on
homogeneous gravitational waves, {\em J. Phys.} A39 (2006)
5189--5226, {\tt arXiv:hep-th/0602036}.

\bibitem{KaprHeis}
M. Kapranov, Heisenberg doubles and derived categories, {\em J.
Algebra} 202, 712--744 (1998), {\tt arXiv:q-alg/9701009}.

\bibitem{Kathotia}
V. Kathotia, Kontsevich's universal formula for deformation
quantization and the Campbell-Baker-Hausdorff formula, {\em Int.
J. Math.}  11:4, 523--551 (2000), {\tt math.QA/9811174}.
%\MR{2002h:53154}

\bibitem{KontsDefPoiss} M. Kontsevich, Deformation quantization
of Poisson manifolds, {\em Lett. Math. Phys.} 66:3, 157--216
(2003)
 %\MR{2005i:53122}

\bibitem{covKappadef}
S. Kre\v si\'c-Juri\'c,  S. Meljanac, M. Stoji\'c, Covariant
realizations of kappa-deformed space, Eur. Phys. J. C 51 (2007),
no. 1, 229–240, {\tt hep-th/0702215}.

\bibitem{Lu}
J-H. Lu, On the Drinfeld double and the Heisenberg double of a
Hopf algebra, {\em Duke Math. J.} 74 (1994) 763–776.

\bibitem{LukNowHeis}
J. Lukierski, A. Nowicki, Heisenberg double description of
$\kappa$-Poincar\'e algebra and $\kappa$-deformed phase space,
{\tt arXiv:q-alg/9792003}.

\bibitem{meljDiff}
S. Meljanac, S. Kre\v{s}i\'c-Juri\'c, Differential structure on
kappa-Minkowski space, and kappa-Poincare algebra, {\em Int. J.
Mod. Physics A} 42:36, 365204-365225 (2009) {\tt arxiv/1004.4647}.

\bibitem{kappaind}
S. Meljanac, A. Samsarov, M. Stoji\'c, K. S. Gupta,
Kappa-Minkowski space-time and the star product realizations, {\em
Eur. Phys.J.} C51 (2007) 229--240, {\tt arXiv:0705.2471}.

\bibitem{scopr} S. Meljanac, Z. \v{S}koda, Leibniz rules for enveloping algebras,
{\tt arXiv:0711.0149}; with newer version at {\tt
http://www.irb.hr/korisnici/zskoda/scopr6.pdf}.

\bibitem{MS} S. Meljanac,
M. Stoji\'{c}, New realizations of Lie algebra kappa-deformed
Euclidean space, {\em Eur.Phys.J.} C47 (2006) 531--539; {\tt
hep-th/0605133}.

\bibitem{exp}
S. Meljanac, D. Svrtan, Z. \v{S}koda, Exponential formulas and Lie
algebra type star products, {\tt arXiv:1006.0478}.

\bibitem{ReshQLieBialg}
N. Reshetikhin, Quantization of Lie bialgebras, {\em Int. Math.
Res. Notices} 7, 143--151 (1992).

\bibitem{Semikhatov}
 A. M. Semikhatov, Heisenberg double addition to the logarithmic
Kazhdan–-Lusztig duality, {\em Lett. Math. Phys.} 92 (2010), no.
1, 81–98, {\tt arXiv:0905.2215}

\bibitem{exterior}
Z. \v{S}koda, Twisted exterior derivative for enveloping algebras,
{\tt arXiv:0806.0978}.

\end{thebibliography}
\end{document}